\documentclass[english]{amsart}
\usepackage{amsmath}
\usepackage{amssymb}
\usepackage{amsthm}
\usepackage{babel}
\usepackage{graphicx}
\newtheorem{lem}{Lemma}
\newtheorem{prop}{Proposition}

\theoremstyle{definition}
\newtheorem{exa}{Example}

\author{Antonio M. Oller Marc\'{e}n}
\title{Counting domino trains}
\address{Departamento de Matem\'{a}ticas, Universidad de Zaragoza\\ C/Pedro Cerbuna 12, 50009 Zaragoza (Spain)}
\email{oller@unizar.es}
\date{}

\def\M{\mathcal{M}}
\def\R{\mathbb{R}}

\begin{document}
\maketitle
\begin{abstract}
In this paper we present a way to count the number of trains that we can construct with a given set of domino pieces. As an application we obtain a new method to compute the total number of eulerian paths in an undirected graph as well as their starting and ending vertices.
\end{abstract}

\section{Introduction}
    Let $G$ be a directed graph with set of vertices $V_G$ and edge set $E_G$. A path of length $m$ on $G$ is a sequence of vertices $v_0,v_1,\dots,v_m$ such that $(v_{i-1},v_i)\in E_G$ for all $1\leq i\leq m$. A path is a cycle if $v_0=v_m$ and $(v_{i-1},v_i)\neq(v_{j-1},v_j)$ for all $1\leq i<j\leq m$. An eulerian path or cycle is a path or cycle of length $|E_G|$. These concepts can be extended to the case when $G$ is an undirected graph in a natural way (see \cite{WES} for these and other elementary concepts about Graph Theory).

    Two eulerian cycles are called equivalent if one is a cyclic permutation of the other. Let $\textrm{Eul}(G)$ denote the number of equivalence classes of eulerian circuits. If $G$ is a directed graph there is a well-known theorem, the so called BEST Theorem (see \cite{TTS,AEB}), which computes $\textrm{Eul}(G)$ but if $G$ is an undirected graph the situation is much more difficult. Nevertheless, in the case of complete graphs interesting results exist. In \cite{McK} for instance, an asymptotic value for $\textrm{Eul}(K_n)$ (as well as the exact number for $n\leq21$) is given, with $K_n$ being the complete graph with $n$ vertices.

    Now let us denote by $\textrm{Eul}_i^j(G)$ the number of eulerian paths in $G$ starting in vertex $v_i$ and ending in vertex $v_j$ (with no equivalence relation taken into account). In this paper we introduce a new approach and present a new method to compute $\textrm{Eul}_i^j(G)$ for any undirected graph $G$. In particular we will count the number of trains that we can construct with a given set of domino pieces, where a train is a chain constructed following the rules of domino and using all the pieces from the given set. Since any graph gives rise to a set of domino pieces where trains correspond to eulerian paths, our method applies to any graph.

    The paper is organized as follows. In the first section we present some elementary definitions and fix the notation. The second section is devoted to prove the main result in the paper, namely to present a method to compute the number of trains constructible from a given set of domino pieces. Finally, in the third section we translate this result to a graph theory setting.

\section{Definitions and notation}
    Let us denote by $\{e_{ij}\ |\ 1\leq i,j\leq n\}$ the set of matrix units in $\M_n(\R)$; i.e., $e_{ij}$ is the square matrix of size $n$ having a 1 in position $(i,j)$ and 0 elsewhere.
    Now, we define $\overline{e}_{ij}=e_{ij}+e_{ji}$ for all $i\neq j$ and $\overline{e}_{ii}=e_{ii}$. Clearly the set $\mathcal{B}=\{\overline{e}_{ij}\ |\ 1\leq i\leq j\leq n\}$ is a basis for the vector space $\mathcal{S}_n(\R)$ of real symmetric matrices of size $n$.
    If we define a new product over $\M_n(\R)$ by $A\bullet B=AB+BA$ our basis multiply in the following way:

    \begin{equation}\label{tabla}\overline{e}_{ij}\bullet\overline{e}_{kl}=\begin{cases}
    \overline{e}_{jl} & \textrm{if $i=k$ and $j\neq l$}\\
    2\overline{e}_{ii}+2\overline{e}_{jj} & \textrm{if $i=k$ and $j=l$ and $i\neq j$}\\
    2\overline{e}_{ii} & \textrm{if $i=j=k=l$}\\
    0 & \textrm{if $\{i,j\}\cap\{k,l\}=\emptyset$}
    \end{cases}\end{equation}

    It is well-known that $\bullet$ is commutative but non-associative. Nevertheless, we may define recursively $A_1\bullet\dots\bullet A_r=(A_1\bullet\dots\bullet A_{r-1})\bullet A_r$. Now, given $m\in\mathbb{N}$ we consider the polynomial $S_m(X_1,\dots,X_m)=\displaystyle{\sum_{\sigma\in S_m}X_{\sigma(1)}\cdots X_{\sigma(m)}}$. With the previous convention and
    choosing $\overline{e}_t=\overline{e}_{i_tj_t} $ for $t=1,\dots,m$ it makes sense to compute $S_m(\overline{e}_1,\dots,\overline{e}_t)=\displaystyle{\sum_{\sigma\in S_m}\overline{e}_{\sigma(1)}\bullet\cdots\bullet\overline{e}_{\sigma(m)}}$. Moreover, there must exist $\{\alpha_{ij}\}\subseteq\R$ such that $S_m(\overline{e}_1,\dots,\overline{e}_m)=\sum\alpha_{ij}\overline{e}_{ij}$.

    In what follows we will denote by $(i,j)$ the domino piece marked with numbers $i$ and $j$. We recall that in the game of domino two pieces can be placed together if they share at least one of their numbers. Now let us suppose that we are given certain set of dominoes that we shall denote by $\mathcal{D}=\{(i_1,j_1),\dots,(i_m,j_m)\}$, we define a train as a sequence $(i_{k_1},j_{k_1})\dots(i_{k_m},j_{k_m})$ admissible by the rules of domino; i.e., such that $j_{k_r}=i_{k_{r+1}}$ for all $1\leq r\leq m-1$.

\section{Counting domino trains}
    Given a set of domino pieces $\mathcal{D}=\{(i_1,j_1),\dots,(i_m,j_m)\}$ we are interested in counting the number of trains that we can construct using all the pieces of $\mathcal{D}$. If an element of $\mathcal{D}$ appears more than once, we will assume that we can distinguish them.

     For any domino piece $(i,j)$ we will identify $(i,j)\leftrightarrow\overline{e}_{ij}$. Clearly, two pieces $(i,j)$ and $(k,l)$ can be placed together following the rules of domino if and only if $\overline{e}_{ij}\bullet\overline{e}_{kl}\neq 0$. Consequently we have the following:

    \begin{lem}\label{lema1}
    A sequence $(i_1,j_1)(i_2,j_2)\dots(i_n,j_n)$ is a train if and only if  $\overline{e}_{i_1j_1}\bullet\overline{e}_{i_2j_2}\bullet\cdots\bullet\overline{e}_{i_nj_n}\neq0$.
    \end{lem}
    \begin{proof}
    By induction on $n$. Cases $n=1,2$ are obvious due to (\ref{tabla}). Now let us suppose that  $(i_1,j_1)(i_2,j_2)\dots(i_n,j_n)$ is a train with $n\geq 3$ and $\overline{e}_{i_1j_1}\bullet\overline{e}_{i_2j_2}\bullet\cdots\bullet\overline{e}_{i_nj_n}=0$. Then, there must exist $2\leq k<n$ such that $\overline{e}_{i_1j_1}\bullet\cdots\bullet\overline{e}_{i_kj_k}=0$ and by our induction hypothesis $(i_1,j_1)\dots (i_k,j_k)$ is not a train which is a contradiction. Conversely if $\overline{e}_{i_1j_1}\bullet\overline{e}_{i_2j_2}\bullet\cdots\bullet\overline{e}_{i_nj_n}\neq0$, then $\overline{e}_{i_1j_1}\bullet\overline{e}_{i_2j_2}\neq0$ and $\overline{e}_{i_2j_2}\bullet\cdots\bullet\overline{e}_{i_nj_n}\neq0$ so, by induction, both $(i_1,j_1)(i_2,j_2)$ and $(i_2,j_2)\dots (i_n,j_n)$ are trains and consequently so is $(i_1,j_1)(i_2,j_2)\dots(i_n,j_n)$.
    \end{proof}

    Now, given the set $\mathcal{D}$ and since the set $\mathcal{B}$ is a basis for $\mathcal{S}_n(\R)$, we have that $S_m(\overline{e}_{i_1j_1},\dots,\overline{e}_{i_mj_m})=\sum\alpha_{ij}\overline{e}_{ij}$. In fact, the following lemma holds.

    \begin{lem}
    In the previous situation, $\alpha_{i_kj_k}\neq0$ if and only if we can construct a train starting with $i_k$ and ending with $j_k$ (or viceversa) using the pieces of the set $\mathcal{D}$.
    \end{lem}
    \begin{proof}
    It is an easy consequence of (\ref{tabla}) and Lemma \ref{lema1}.
    \end{proof}

    Observe that, given a train $(i_1,j_1)(i_2,j_2)(i_3,j_3)$ we could have placed its pieces in different ways according to the order. Namely we could have firstly placed the piece $(i_1,j_1)$ then $(i_2,j_2)$ on its right and finally the piece $(i_3,j_3)$ on the right of the latter. We could also have started by piece $(i_2,j_2)$, then $(i_1,j_1)$ on its left and finally $(i_3,j_3)$ on the right of the former. If we denote by $c(n)$ the number of different ways we can construct a train of length $n$ we have.

    \begin{lem}
    $c(n)=2^{n-1}$
    \end{lem}
    \begin{proof}
    Clearly we can finish a given train either by placing the first piece (the one in the left) or the last one (the one in the right). This implies that $c(n)=2c(n-1)$ and since $c(1)=1$ the result follows.
    \end{proof}

    With these technical lemmas in mind, we are in condition to prove the main result of this paper.

    \begin{prop}\label{main}
    Given a set of domino pieces $\mathcal{D}=\{(i_1,j_1),\dots,(i_m,j_m)\}$, the number $\displaystyle{\frac{\alpha_{i_kj_k}}{2^{m-1}}}$ is exactly the number of domino trains that we can construct starting with $i_k$ and ending with $j_k$ (or viceversa) using the pieces from the set $\mathcal{D}$.
    \end{prop}
    \begin{proof}
    Again, cases $m=1,2$ are obvious and can be verified by direct computation using (\ref{tabla}). Now if $m\geq 3$ it is enough to observe that $$\displaystyle{S_m(\overline{e}_{i_1,j_1},\dots,\overline{e}_{i_m,j_m})=\sum_{k=1}^{m}S_{m-1}(\overline{e}_{i_1j_1},\dots,\widehat{\overline{e}_{i_kj_k}},\dots,\overline{e}_{i_mj_m})\bullet\overline{e}_{i_kj_k}}$$
    and proceed by induction on $m$.
    \end{proof}

    \begin{exa}
    Given the set $\mathcal{D}=\{(1,1),(1,2),(1,3),(2,2),(2,3),(3,3)\}$ we have that $S_6(\overline{e}_{11},\overline{e}_{12},\overline{e}_{13},\overline{e}_{22},\overline{e}_{23},\overline{e}_{33})=2^5\left(\displaystyle{\sum_{i=1}^34\overline{e}_{ii}}\right)$.
    So, we can construct 4 trains starting and ending with 1, and the same number with 2 and 3.
    \end{exa}

    \begin{exa}
    Given the set $\mathcal{D}=\{(1,2),(1,3),(2,3),(2,3),(2,4),(3,4)\}$ we have that $S_6(\overline{e}_{12},\overline{e}_{13},\overline{e}_{23},\overline{e}_{22},\overline{e}_{24},\overline{e}_{34})=2^5\left(12\overline{e}_{11}+24\overline{e}_{22}+24\overline{e}_{33}+12\overline{e}_{44}\right)$.
    So, we can construct 12 trains starting and ending with 1 (the same number with 4), and 24 trains with 2 (the same number with 3).
    \end{exa}

\section{Counting eulerian paths in graphs}

    Given a set of dominoes $\mathcal{D}=\{(i_1,j_1),\dots(i_m,j_m)\}$ we can construct a graph in the following straightforward way: we draw a vertex for each of the numbers appearing ir our set $\mathcal{D}$ and we join two of such vertices if and only if the corresponding piece lies in our set. Clearly we can proceed backwards and obtain a set of dominoes from every graph. Moreover, domino trains correspond exactly with eulerian paths or cycles. Thus, our previous construction allows us to count them.

     Also observe that given a graph $G$ if we denote by $\mathcal{D}_{G}$ the set of dominoes coming from the previous construction, then $\displaystyle{\sum_{(i,j)\in\mathcal{D}_G}\overline{e}_{ij}=A_ G}$ where $A_G$ is the adjacency matrix of $G$. Conversely, given a graph $G$ with adjacency matrix $A_G$ we can write $A_G=\displaystyle{\sum_{k=1}^m\overline{e}_{i_kj_k}}$ with $m=|E_G|$ and we know that $S_m(\overline{e}_{i_1j_1},\dots,\overline{e}_{i_mj_m})=\sum\alpha_{ij}\overline{e}_{ij}$. The following result translates Proposition \ref{main} to this setting.

     \begin{prop}
     In the previous situation and with the notation of Section 1, $\textrm{Eul}_i^j(G)=\displaystyle{\frac{\alpha_{ij}}{2^{m-1}}}$.
     \end{prop}

    \begin{exa}
    Given the graph in the figure
    \begin{figure}[h]
    \begin{center}
    \includegraphics[width=5 cm]{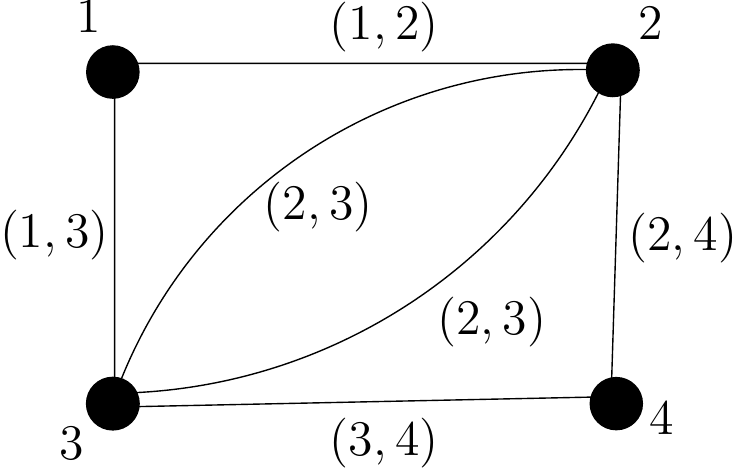}
    \end{center}
    \end{figure}
     it is clear that the associated set of domino pieces is $\mathcal{D}_G=\{(1,2),(1,3),(2,3),(2,3)(2,4),(3,4)\}$ and according to Example 2 above, the total number of eulerian cycles of $G$ is 72, 12 starting and ending in vertex 1, 24 in vertex 2, 24 in vertex 3 and 12 in vertex 4.
    \end{exa}

    \begin{exa}
    Consider the graph $G$ given in the figure below.
    \begin{figure}[h]
    \begin{center}
    \includegraphics[width=5 cm]{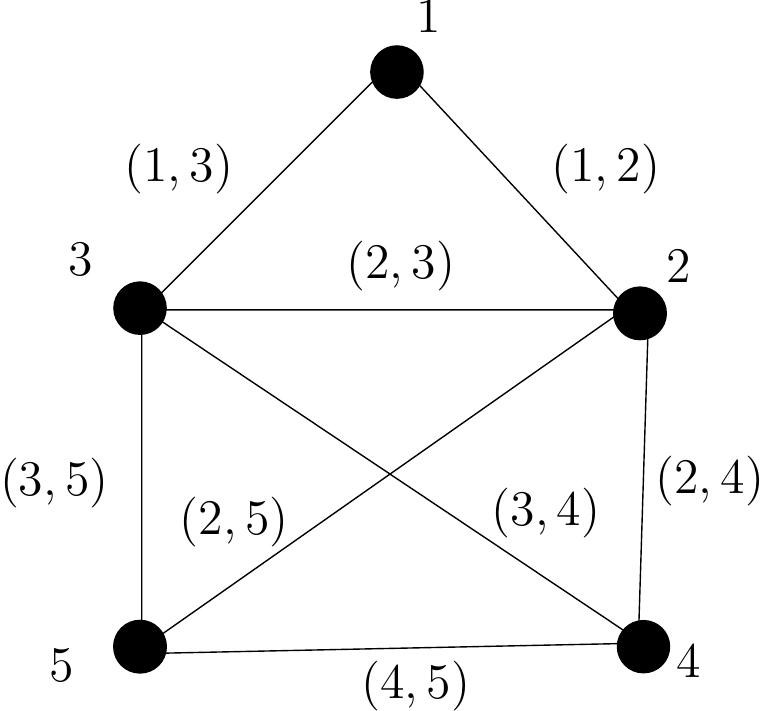}
    \end{center}
    \end{figure}

    In this case $S_8(\overline{e}_{12},\overline{e}_{13},\overline{e}_{23},\overline{e}_{24},\overline{e}_{25},\overline{e}_{34},\overline{e}_{35},\overline{e}_{45})=2^7\cdot 44\overline{e}_{45}$ and thus $G$ has 44 eulerian paths starting in vertex 4 and ending in vertex 5 (and vice versa).
    \end{exa}





  \section*{Acknowledgments}

    The author thanks the Spanish project MTM2007-67884-C04-02
    for the partial financial support.

\bibliography{./ref}

\begin{thebibliography}{1}

\bibitem{McK}
Brendan~D. McKay and Robert~W. Robinson.
\newblock Asymptotic enumeration of {E}ulerian circuits in the complete graph.
\newblock {\em Combin. Probab. Comput.}, 7(4):437--449, 1998.

\bibitem{TTS}
W.~T. Tutte and C.~A.~B. Smith.
\newblock On {U}nicursal {P}aths in a {N}etwork of {D}egree 4.
\newblock {\em Amer. Math. Monthly}, 48(4):233--237, 1941.

\bibitem{AEB}
T.~van Aardenne-Ehrenfest and N.~G. de~Bruijn.
\newblock Circuits and trees in oriented linear graphs.
\newblock {\em Simon Stevin}, 28:203--217, 1951.

\bibitem{WES}
Douglas~B. West.
\newblock {\em Introduction to graph theory}.
\newblock Prentice Hall Inc., Upper Saddle River, NJ, 1996.

\end{thebibliography}
\bibliographystyle{plain}

\end{document}